\documentclass[12pt]{amsart}
\usepackage{amsthm}
\usepackage{amsmath, mathrsfs}
\usepackage{amssymb}
\usepackage{hyperref}
\usepackage{enumerate}
\usepackage{latexsym,array}
\usepackage{amsfonts}
\usepackage{shadow}

\newtheorem{Pa}{Paper}[section]
\newtheorem{Tm}[Pa]{{\bf Theorem}}
\newtheorem{La}[Pa]{{\bf Lemma}}

\newtheorem{Pn}[Pa]{{\bf Proposition}}

\newcommand{\va}{\varphi}
\date{}
\author[D. Alpay]{Daniel Alpay}
\author[P. Jorgensen]{Palle Jorgensen}
\address{(DA) Department of Mathematics
\newline
Ben Gurion University of the Negev \newline P.O.B. 653,
\newline
Be'er Sheva 84105, \newline ISRAEL}
\email{dany@math.bgu.ac.il}
\address{(PJ)
Department of Mathematics\newline 14 MLH \newline The University
of Iowa Iowa City,\newline IA 52242-1419 USA}
\email{palle-jorgensen@uiowa.edu}
\thanks{Daniel Alpay wishes to thank the
Earl Katz family for endowing the chair which supported his
research. The first three named authors wish to thank the
Binational Science Foundation Grant number 2010117.}

\author[I. Lewkowicz]{Izchak Lewkowicz}
\address{(IL) Department of Electrical Engineering
\newline
Ben Gurion University of the Negev \newline P.O.B. 653,
\newline
Be'er Sheva 84105, \newline ISRAEL}
\email{izchak@ee.bgu.ac.il}

\author[I. Martziano]{Itzik Martziano}
\address{(IM) Ben Gurion University of the Negev \newline P.O.B. 653,
\newline
Be'er Sheva 84105, \newline ISRAEL}
\email{martzian@bgu.ac.il}
\title[Infinite product and iterations of functions]
{Infinite product representations for kernels and iterations of functions}
\begin{document}
 \dedicatory{Dedicated to the memory of Bill
Arveson}

\begin{abstract}
 We study infinite products of reproducing kernels with view to their use in dynamics
(of iterated function systems), in harmonic analysis, and in
stochastic processes. On the way, we construct a new family of
representations of the Cuntz relations. Then, using these
representations we associate a fixed filled Julia set with a
Hilbert space. This is based on analysis and conformal geometry
of a fixed rational mapping $R$ in one complex variable, and its
iterations.
\end{abstract}

\maketitle
\tableofcontents
\section{Introduction}
\setcounter{equation}{0}
The purpose of this paper is twofold, one is to offer a general
framework for an harmonic analysis of reproducing kernel Hilbert spaces,
and the other is to explore its applications. On the first point, we offer a
general tool for analysis of positive definite kernels, and their associated
reproducing kernel Hilbert spaces (RKHS). Our analysis is based on a construction
of families of representations of a system of operators. We use these
representations in order to construct explicit factorizations of the
kernels at hand; and,  as a result we get computable decompositions of the
associated reproducing kernel Hilbert spaces. Our starting point is a
symbolic set with relations, often called the Cuntz relations (CR) after
the $C^*$-algebra they generate, see \cite{Cun77}. The CRs are indexed by the
number  $n$  of symbols in a generating system; for every $n$ , including
possibly $n =  \infty$, we have a Cuntz algebra  $\mathcal O_n$, see  \cite{Cun77}.\\

Our present results were announced in \cite{aj_pnas}.  While representations of the Cuntz algebras
$\mathcal O_n$, denoted ${\rm Rep}\,(\mathcal O_n, \mathcal H)$
 for a given separable Hilbert space $\mathcal H$, have a long history in operator algebras and their applications,
our present use of them in an analysis of dynamics of complex substitution-systems is novel. The use of
${\rm Rep}\,(\mathcal O_n, \mathcal H)$ in operator algebras and physics was pioneered by Arveson \cite{MR987590};
see also \cite{MR1469149}, and the references there. Fix $n$, then, up to a natural action of the group  $U_n$,
the space ${\rm Rep}\,(\mathcal O_n, \mathcal H)$  “is”  ${\rm End}_n(\mathbf B(\mathcal H))$
where $\mathbf B(\mathcal H))$ denotes the set of all bounded operators on $\mathcal H$ and
${\rm End}_n$  denotes the set of all endomorphisms of $\mathbf B(\mathcal H)$  of Powers index $n$.\\

{\bf 1.1 Representations.}
Our present use of the CRs centers on an analysis of representations of
$\mathcal O_n$;  so our focus is on the representations as opposed to $C^*$-algebraic
questions.  By a representations of the CRs we mean a multivariable operator system
which satisfies the formal Cuntz relations. The study of representations of the Cuntz
relations is subtle for a number of reasons. For example,  it is known (see e.g.,
\cite{BrJo02a})  that, in general, for  $n$ fixed, that the variety of all
representations of $\mathcal O_n$  (up to unitary equivalence) is very large;
to be precise, in fact the set of equivalence classes of irreducible representations
of $\mathcal O_n$  is known to not admit a measurable parameterization with any
Borel set.
   Nonetheless, for diverse purposes, there are interesting infinite
   dimensional families of representations serving diverse purposes in
   harmonic analysis and in applications. Our focus here is a set of
   applications of representations to an harmonic analysis of kernel Hilbert
   spaces. With this, we then arrive at decompositions which are of relevance
  in the study of sub-band filters in signal processing, in an harmonic analysis of de Branges spaces,
and of use in building  $\mathbf L_2$-spaces on fractals, and in an analysis on Hilbert spaces built on Julia sets
arising in conformal geometry from iteration of fixed rational mappings.\\

{\bf 1.2 Kernels.}
The decomposition of a positive definite kernel into a sum of
such kernels is not unusual, but the decomposition of a positive
definite kernel into a product of positive definite kernels seems much less common. In
the work \cite{ajlm1}, such a multiplicative decomposition was
used in the setting of the Hardy space of the open unit disk, and
connections with the Cuntz relations were pointed out. In the
present work, we study positive kernels of the form
\[
K(z,w)=\prod_{n=0}^N k(R_n(z),R_n(w)),
\]
and
\[
K(z,w)=\prod_{n=0}^\infty k(R_n(z),R_n(w)).
\]
In these expressions, $k(z,w)$ is a positive definite function on a
set $\mathcal U$, and $R$ is a map from $\mathcal U$ into itself.
We denote by $R_0(z)=z$ and
\begin{equation}
\label{eqrn}
R_n=\underbrace{R\circ R\circ\cdots\circ
R}_{\mbox{\rm $n$ iterations}},
\end{equation}
and appropriate assumptions are made to insure convergence.\\



{\bf 1.3 Local analysis.}
Our paper is at the crossroads of complex dynamics and
representation theory (representations of certain non-abelian
algebras). To aid readers from both areas, we have included some
preliminaries on the two subjects. Beginning with complex
dynamics, recall that Julia sets (and their compliments Fatou
sets  $F(f)$ ) serve to chart a number of geometric patterns of
points under iterated substitution. Specifically, the two
complementary sets are defined from a given (and fixed) function
$f$;  often  $f$ is a rational function defined on the Riemann
sphere. Intuitively, the Fatou set of  $f$ is the set of points
such that all nearby values under  $f$ behave similarly under
repeated iteration (self-substitution): $F(f)$ is open, and  the
iterations form a normal family in $F(f)$ in the sense of Montel.
By contrast, the points in the compliment, the Julia set  $J(f)$
have the property that an arbitrarily small perturbation can
cause drastic changes in the sequence of iterated function
values, i.e., points in $J(f)$ represents chaotic behavior. The
bounded connected component of  $F(f)$ is called the filled Julia
set, and  $J(f)$ is its boundary. In these dynamical terms, the
intuition is that points in the Fatou set are 'regular', and the
Julia set represents  'chaos'. In terms of geometry of repeated
iteration, the Julia set falls in the realm between
deterministic and stochastic. As a result, the Julia set is
aperiodic, and so does not lend itself to standard tools of
harmonic analysis. Our aim is to zoom in on local properties of
points in the filled Julia set.  It is our aim to show, for
certain examples (some cases when  $f$ is a polynomial),  that
there is a local version of a harmonic analysis which works
well in a neighborhood of attracting cycles in the filled Julia
set, or in an open basin of attraction for an attractive fixed
point of  $f$. Due to a theorem of Brolin \cite{MR0194595}, every given  $f$ has an
invariant maximal-entropy measure $\mu$ (depending on $f$) , and
the support of $\mu$  is $J(f)$. Intuitively, in our analysis ,
$\mu$ plays a role analogous to that played by Haar measure in the
harmonic analysis of compact groups.  A key tool in our analysis
is the design of representations of a system of isometries in the
Hilbert space  $H:= \mathbf L_2(\mu)$  defined from  $\mu = \mu(f)$. Such
a system of  $N$  isometries was defined axiomatically by J. Cuntz
(for entirely different purposes), and when $N$ is fixed we speak
of the Cuntz algebra  $O_N$.\\

In general, the theory of the
representations of $O_N$ is difficult. Nonetheless, for our
present purpose, we show that there is a distinguished subclass of
representations, important for our local harmonic analysis in the
filled Julia set. To help understand this, note that from every
representation of $O_N$ acting on a fixed Hilbert space  $H$ , one
naturally obtains an associated  N-nary system of mutually
orthogonal subspaces in $H$. This is a Hilbert space theoretic
tree-like symbolic representation of local features of the
dynamics. We will construct our representations such that the
corresponding  $N$-nary system in $H= \mathbf L_2(\mu)$  corresponds to
the geometry of iterated backward
     substitutions for a fixed polynomial  $f$.

\section{Preliminaries on reproducing kernel Hilbert spaces}
\setcounter{equation}{0}
In this section we review some facts on positive definite
functions and on operators between reproducing kernel Hilbert
spaces. See for instance \cite{aron,saitoh}. The knowledgeable reader can proceed directly to Section
\ref{sec3}. We divide the section into two subsections, devoted to
definitions, and $\mathcal L(\varphi)$ de Branges spaces.
\subsection{Generalities}
First recall that a function $K(z,w)$
defined on a set $\Omega$ is positive definite if for every
choice of $M\in\mathbb N$, of points $w_1,\ldots, w_M\in\Omega$
and $c_1,\ldots, c_M\in \mathbb C$, it holds that
\[
\sum_{k,\ell=1}^M c_\ell^*K(w_\ell, w_k)c_k\ge 0.
\]
Associated to a positive definite function is a unique Hilbert
space $\mathcal H(K)$ of functions on $\Omega$, with the
following two properties: For every $w\in\Omega$, the function
$K_w\,:\, z\mapsto K(z,w)$ belongs to $\mathcal H(K)$, and for
every $f\in\mathcal H(K)$,
\[
\langle f,K_w\rangle_{\mathcal H(K)}=f(w).
\]
The function $K$ is called the reproducing kernel of the space
and it can be computed by the formula
\begin{equation}
\label{eqrk}
K(z,w)=\sum_{j\in J}e_j(z)e_j(w)^*,
\end{equation}
where $(e_j)_{j\in J}$ is any orthonormal basis of $\mathcal
H(K)$. Note that $J$ need not be countable since the space
$\mathcal H(K)$ need not be separable.\\

In the case of the Hardy space $\mathbf
H_2(\mathbb C_+)$ of the open right-half-plane $\mathbb C_+$, an
orthonormal basis for $\mathbf H_2(\mathbb C_+)$ is given by the
functions
\[
t_n(z)=\frac{1}{\sqrt{\pi}}\frac{1}{z+1}\left(\frac{z-1}{z+1}\right)^{n},\quad
n=0,1,\ldots
\]
and we have
\[
\frac{1}{2\pi(z+w^*)}=\sum_{n=0}^\infty t_n(z)t_n(w)^*.
\]

\begin{Pn}
Let $K$ be positive definite on $\Omega$ and let $\varphi$ denote
a function from $\Omega$ into itself, and let $e$ be a function
from $\Omega$ into $\mathbb C$. The operator
\[
Sf(z)=e(z)f(\va(z))
\]
is bounded from $\mathcal H(K)$ into itself if and only if the
function
\begin{equation}
\label{kephi}
K(z,w)-e(z)e(w)^*K(\va(z),\va(w))
\end{equation}
is positive on $\Omega$. When this condition holds, the adjoint
operator is given by the formula
\begin{equation}
\label{mutl}
S^*K_w=e(w)^*K_{\varphi(w)}.
\end{equation}
\end{Pn}

While the Proposition can be found in the literature, we include below a sketch of the main idea involved as
it serves to unify several themes coming later in a variety of seemingly different context;
harmonic analysis, representation theory, the study of Julia sets, to mention a few.
The details below further serve to introduce terminology to be used later. To prove the proposition consider the linear
relation in $\mathcal H(K)\times\mathcal H(K)$ spanned by the pairs
\[
(K(\cdot, w),e(w)^*K(\cdot, \varphi(w)),\quad w\in\Omega.
\]
It is densely defined. Moreover it is contractive thanks to \eqref{kephi}. Therefore it extends to the graph
of an everywhere defined contraction, say $T$. For $f\in\mathcal H(K)$ and $w\in \Omega$ we have
\[
\begin{split}
\langle T^*f,K(\cdot, w)\rangle_{\mathcal H(K)}&=\langle f,T(K(\cdot, w))\rangle_{\mathcal H(K)}\\
&=\langle f, e(w)^*K(\cdot, \varphi(w))\rangle_{\mathcal H(K)}\\
&=e(w)f(\varphi(w))\\
&=(Sf)(w).
\end{split}
\]

\subsection{$\mathcal L(\varphi)$ spaces}
Recall that a $\mathbb
C^{n\times n}$-valued function $\varphi$ is analytic in $\mathbb
C_+$ and such that ${\rm Re}\,\varphi(z)\ge 0$ for $z\in\mathbb C_+$
if and only if it is of the form
\begin{equation}
\label{weizmann} \varphi(z)=a+bz-i\int_{\mathbb
R}d\mu(t)\left\{\frac{1}{t-iz}-\frac{t}{t^2+1}\right\},
\end{equation}
where $a\in \mathbb C^{n\times n}$ is such that $a+a^*=0$,
$b\in\mathbb C^{n\times n}$ is non negative, and where $d\mu$ is
a $\mathbb C^{n\times n}$-valued positive measure subject to
\[
\int_{\mathbb R}\frac{d\mu(t)}{t^2+1}<\infty
\]
This expression allows to extend $\varphi$ to the open left half
plane (the extension will not be continuous across the imaginary
axis in general). When $\varphi$ is extended in such a way, the
kernel
\begin{equation}
\label{kernphi}
\frac{\varphi(z)+\varphi(w)^*}{z+w^*}=b+\int_{\mathbb
R}\frac{d\mu(t)}{(t-iz)(t-iw)^*}
\end{equation}
is positive definite in $\mathbb C\setminus i\mathbb R$. The
associated reproducing kernel Hilbert space $\mathscr L(\varphi)$ was
characterized and studied by de Branges. It consists of functions
of the form
\[
F(z)=b\xi+\int_{\mathbb R}\frac{d\mu(t)f(t)}{t-iz},
\]
with $\xi\in\mathbb C^n$ and $f\in\mathbf L_2(d\mu)$, with norm
\[
\|F\|^2=\xi^*b\xi+\|f\|_{\mu}^2.
\]
See \cite{MR0229011}, \cite{ad1}. The space is finite dimensional if and
only if $\mu$ is a jump measure with a finite number of jumps.
This happens if and only if the function $\phi$ is rational and
satisfies
\[
\varphi(z)=-\varphi(-z^*)^*,\quad \forall z\in\mathbb C,\,\mbox{\rm
which is not a pole of $\phi$.}
\]
We now consider $n=1$ and consider a rational function $\va$ such
that $\mathscr L(\va)$ is finite dimensional. Thus
\[
\varphi(z)=a+bz+\sum_{j=1}^N\frac{m_j}{t_j-iz}
\]
where the $t_j\in\mathbb R$ and the $m_j>0$.

\begin{Pn}
Let $a=b=0$. Let
\[
e_j(z)=\frac{\sqrt{m_j}}{(t_j-iz)},\quad j=1,\ldots N.
\]
Then,
\begin{equation}
\label{phi1}
\frac{\va(z)+\va(w)^*}{z+w^*}=\sum_{n=1}^Ne_n(z)e_n(w)^*,
\end{equation}
and
\begin{equation}
\label{phi2}
\langle e_j,e_k\rangle_{\mathscr
L(\phi)}=\delta_{jk}.
\end{equation}
\end{Pn}

{\bf Proof:} Equation \eqref{phi1} is a special case of
\eqref{kernphi} when $\mu$ is a jump measure with a finite number
of jumps. Formula \eqref{phi2} is a special case of \eqref{eqrk}
since the $e_j$ are linearly independent.
\mbox{}\qed\mbox{}\\

\section{A general setting}
\label{sec3}
We now assume that the set $\mathcal U$ in the introduction is a topological space and
obtain infinite product representations of certain positive definite kernels.
We start from a function $k(z,w)$ positive definite in $\mathcal U$
and denote by $(e_j)_{j\in J}$ denote an orthonormal basis of the
reproducing kernel Hilbert space $\mathcal H(k)$ with reproducing
kernel $k(z,w)$ (with $z,w\in\mathcal U$). Thus,
\begin{equation}
\label{formule_k}
k(z,w)=\sum_{j\in J}e_j(z)(e_j(w))^*, \quad
z,w\in\mathcal U.
\end{equation}
At this stage, $\mathcal H(k)$ need
not be separable, and thus the index set $J$ need not be countable.
We assume that
$\mathcal U$ satisfies
\begin{equation}
K(z,w)=k(z,w)K(R(z),R(w)),\quad \forall z,w\in\mathcal U.
\label{2}
\end{equation}
where the function $K(z,w)$ is not {\it a priori} not positive definite in $\mathcal U$.

\begin{Pn}
Let $k(z,w)$ positive definite in $\mathcal U$, and
assume that \eqref{2} is in force for some function $K(z,w)$, which is continuous on $\mathcal U$, and not identically
equal to $0$. Assume that there exists a point
$\ell\in\mathcal U$ such that
\begin{equation}
\label{fixpoint}
\lim_{n\rightarrow\infty} R_{n}(z)=\ell,\quad\forall z\in\mathcal U.
\end{equation}
Then if $K(\ell,\ell)>0$, the function $K(z,w)$ is positive definite in $\mathcal U$ and
\[
K(z,w)=\left(\prod_{n=0}^\infty\left(\sum_{i\in I}e_i(R_n
(z))e_i(R_n(w))^*\right)\right)K(\ell ,\ell).
\]
\end{Pn}

{\bf Proof:} Let $N\in\mathbb N_0$. It holds that
\begin{equation}
\label{prod1}
K(z,w)=\left(\prod_{n=0}^Nk(R_n(z),R_n(w))\right)K(R_{N+1}(z),R_{N+1}(w)),
\end{equation}
for $z,w\in\mathcal U$. We note that $R(\ell)=\ell$.
The hypothesis imply that
\[
\lim_{N\rightarrow\infty}K(R_{N+1}(z),R_{N+1}(w))=K(\ell,\ell)>0,
\]
and hence the infinite product
$\prod_{n=0}^\infty k(R_n(z),R_n(w))$ converges for $z,w\in\mathcal U$ and is equal to $\frac{K(z,w)}{K(\ell,\ell)}$.
\mbox{}\qed\mbox{}\\

In the preceding proposition, one does not assume that $k(z,w)-1$ is positive definite in $\mathcal U$.
As examples, we mention the works \cite{JP98b}, \cite{}, where
positive definite kernels of the form
\[
\prod_{n=0}^\infty \cos\frac{(t-s)}{4^n}
\]
are introduced. The main point of the example is to illustrate that our method works in examples that have more gaps
than are usually involved in the standard theory for inverse iteration of branches of some fixed polynomial, or rational
function $R(z)$. Starting with $R$ of degree $N$, say, it is natural to create an IFS corresponding to a choice of $N$
branches of inverse for $R$. In more detail, let $R$ be fixed, and let the degree of $R$ be $N$.  Then the Riemann
surface for $R$ has $N$ sheets, and $R$ will be onto with a system of $N$  functions, serving as branches of inverse for
$R$. One then iterates such an  $N$-nary system of inverses (see \cite{DuJo06a}.) Even choosing for $R$
just the monomial  $R(z):=z^4$  leads to IFSs with gaps of interest in harmonic analysis of lacunary Fourier expansions.
Recall that in the gap-examples, such as $R(z)=z^4$,  initially there are four distinct functions as inverse for
$R$, but one may select only two of them for an IFS. The result is a fractal with gaps, and dimension $1/2$.
It can be represented as a Cantor set $J(4,2)$, here realized as a subset of the circle (= the Julia set for $R$).
In this case, the Brolin measure  $\mu$  (see \cite{MR0194595}) coincides with the IFS measure (of dimension $1/2$ )
corresponding to the choice
of two branches of inverse for $R(z)=z^4$. The support of  $\mu$  coincides with Cantor set  $J(4,2)$. We get a Hilbert
space of lacunary power series with $\mathbf L_2(\mu)$-boundary values supported on the Cantor set  $J(4,2)$. We also
refer to \cite{MR2097020,MR2063755,MR2239365,MR2240292,MR2270247,MR2277210} for related works on IFS and CR.

\subsection{A representation of the Cuntz algebra $\mathcal O_N$}

In this subsection we construct representations of the Cuntz
relations. See \cite{MR2277210} for more on these relations.

\begin{Tm}
Let $\mathcal H(K)$ denote the reproducing kernel
of functions defined on $\mathcal U$ and with reproducing kernel
$K(z,w)$.
The operators
\[
(S_jf)(z)=e_j(z)f(R(z)),\quad j\in J,
\]
are continuous from $\mathcal H(K)$ into itself and satisfy
\[
\sum_{j\in J} S_jS_j^*=I_{\mathcal H(K)}.
\]
\end{Tm}

{\bf Proof:} Let $w\in\mathcal U$. From formula \eqref{mutl} we have
\[
S_j^*K_w=(e_j(w))^*K_{R(w)},
\]
and so
\[
(S_jS_j^*K_w)(z)=e_j(z)(e_j(w))^*K(R(z),R(w)).
\]
It follows that
\[
\begin{split}
(\sum_{j\in J} S_jS_j^*K_w)(z)&=\left(\sum_{j\in
J}e_j(z)(e_j(w))^*\right)K(R(z),R(w))\\
&=k(z,w)K(R(z),R(w))\quad (\mbox{\rm using \eqref{formule_k}})\\
&=K(z,w)\quad (\mbox{\rm using \eqref{bastille}}).
\end{split}
\]
\mbox{}\qed\mbox{}\\

Two cases of interest, which will be elaborated upon in the
following section, correspond to $\mathcal U=\mathbb D$, $\varphi(z)=z^2$ and
$k(z,w)=1+zw^*$ and $k(z,w)=(1+zw^*)^2$ respectively. In the
first case, $K(z,w)=\frac{1}{1-zw^*}$ and $\mathcal H(K)$ is
equal to $\mathbf H_2(\mathbb D)$, the Hardy space of the open
unit disk. In the second case, $K(z,w)=\frac{1}{(1-zw^*)^2}$ and
$\mathcal H(K)$ is equal to $\mathbf B_2(\mathbb D)$, the Bergman
space of the open unit disk. In the first case the $S_j$ satisfy
the Cuntz relations  while they do not satisfy these relations in
the second case.

\section{Harmonic analysis of kernels from Blaschke products and Bergmann space}
\setcounter{equation}{0}
\subsection{$\mathcal H(b)$ spaces}
We set
\begin{equation}
\label{eqbl}
b(z)=\prod_{i=1}^N\frac{z-w_i}{1-zw_i^*}
\end{equation}
to be a finite Blaschke product of the open unit disk $\mathbb D$.
Writing
\[
\frac{1}{1-zw^*}=\frac{1}{1-b(z)b(w)^*}\frac{1-b(z)b(w)^*}{1-zw^*},
\]
one obtains a multiplicative decomposition of the Cauchy kernel.
Setting $e_1,\dots, e_N$ to be an orthonormal basis
of $\mathbf H_2\ominus b\mathbf H_2$ we have
\begin{equation}
\label{multi}
\frac{1}{1-zw^*}=(\sum_{i=1}^N e_i(z)e_i(w)^*)\frac{1}{1-b(z)b(w)^*},
\end{equation}
and so, for every $M\in \mathbb N$,
\begin{equation}
\label{multi11}
\frac{1}{1-zw^*}=\left(\prod_{n=0}^M(\sum_{i=1}^N
e_i (b^{\circ n}(z))e_i(b^{\circ
n}(w)^*))\right)\frac{1}{1-b^{\circ (M+1)}(z)(b^{\circ (M+1)}(w))^*}.
\end{equation}
Assume that
\[
w_1=w_2=0
\]
in \eqref{eqbl}. Then,
\begin{equation}
\label{limitzero}
\lim_{M\rightarrow\infty} b^{\circ M}(z)=0,\quad
\forall z\in\mathbb D,
\end{equation}
and we obtain the infinite product representation
\begin{equation}
\label{multi1}
\frac{1}{1-zw^*}=\prod_{n=0}^\infty(\sum_{i=1}^N
e_i (b^{\circ n}(z))e_i(b^{\circ n}(w)^*)),
\end{equation}
where we have denoted
\[
b^{\circ n}(z)=\begin{cases} z,\quad\hspace{2.5cm} {\rm if}\quad n=0,\\
\underbrace{(b\circ b\circ \cdots\circ b)(z)}_{n\,\,{\rm
times}},\quad\hspace{-3mm}{\rm if}\quad n=1,2,\ldots
\end{cases}
\]
From \eqref{multi1} one obtains mutiplicative decompositions for the
kernels $\frac{1}{(1-zw^*)^t}$, $t=2,3,\ldots$.\\


Furthermore, \eqref{multi1} implies the orthogonal decomposition
\[
\mathbf H_2=\oplus_{i=1}^N e_i(\mathbf H_2\ominus b\mathbf H_2)
\]
and the maps $S_if(z)=e_i(z)f(b(z))$ are bounded from
$\mathbf H_2$ into itself and satisfy the Cuntz relations.

\subsection{The Bergmann space}
Let $b$ be the Blaschke product of degree $N$ defined in
\eqref{eqbl}. In the case of the Bergmann space we have
\[
\frac{1}{(1-zw^*)^2}=\frac{1}{(1-b(z)b(w)^*)^2}\frac{(1-b(z)b(w)^*)^2}{(1-zw^*)^2}.
\]
Both the kernels
\[
\frac{1}{(1-b(z)b(w)^*)^2}\quad{\rm and}\quad
\frac{(1-b(z)b(w)^*)^2}{(1-zw^*)^2}
\]
are positive definite in $\mathbb D$. Furthermore, with
$e_i,i=1,\ldots, N,$ being an orthonormal basis of $\mathcal H(b)$
we have
\[
\frac{1}{(1-zw^*)^2}=\sum_{i,j=1}^N
e_i(z)e_j(z)e_i(w)^*e_j(w)^*\frac{1}{(1-b(z)b(w)^*)^2}
\]
which leads to the decomposition
\[
\mathcal B=\sum_{i,j=1}^N e_ie_j \mathcal B(b).
\]
This decomposition will not be orthogonal in general.\\

The case $b(z)=z^2$ is of special interest. Then,
\[
\begin{split}
\frac{1}{(1-zw^*)^2}&=\frac{(1-b(z)b(w)^*)^2}{(1-zw^*)^2}\frac{1}{(1-b(z)b(w)^*)^2}\\
&=(1+2zw^*+z^2(w^*)^2)K(b(z),b(w)),
\end{split}
\]
and we obtain the multiplicative representation of the Bergmann
kernel
\[
\frac{1}{(1-zw^*)^2}=\prod_{n=0}^\infty
(1+2z^{2^n}(w^*)^{2^n}+z^{2^{n+1}}(w^*)^{2^{n+1}})
\]

\subsection{Functions with real positive part}
 Furthermore,
\begin{equation}
\label{paris}
\begin{split}
\frac{I_n}{z+w^*}&=\frac{1}{\va(z)+\va(w)^*}\cdot\frac{\va(z)+\va(w)^*}{z+w^*}I_n\\
&=\sum_{n=1}^Ne_n(z)\frac{I_n}{\va(z)+\va(w)^*}e_n(w)^*\\
&=\sum_{n=1}^N\sum_{m=0}^\infty
e_n(z)t_m(\va(z))t_m(\va(w))^*e_n(w)^*.
\end{split}
\end{equation}
Each of the term
\[
e_n(z)t_m(\va(z))t_m(\va(w))^*e_n(w)^*
\]
is a positive definite function, of rank $1$. The associated
one-dimensional reproducing kernel Hilbert space is spanned by
the function
\[
z\mapsto e_n(z)t_m(\va(z)).
\]
These spaces do not intersect since for $(n_1,m_1)\not=(n_2,m_2)$
\begin{equation}
\label{qwerty}
ae_{n_1}(z)t_{m_1}(z)+be_{n_2}(z)t_{m_2}(\va(z))\equiv
0\Longrightarrow a=b=0.
\end{equation}

By $(\mathbf H_2(\mathbb C_+))(\va)$ the reproducing kernel
Hilbert space with reproducing kernel $\frac{1}{\va(z)+\va(w)^*}$.

\begin{Pn}
$f\in\mathbf H_2(\mathbb C_+))(\va)$ if and only if it can be
written as
\begin{equation}
\label{paris_passy}
 f(z)=h(\va(z)),\quad h\in\mathbf H_2(\mathbb
C_r)
\end{equation}
with norm
\begin{equation}
\|f\|=\|h\|. \label{paris_trocadero}
\end{equation}
\end{Pn}

{\bf Proof:} We first note that \eqref{paris_trocadero} indeed
defines a quadratic norm on the linear span of functions of the
form \eqref{paris_passy}, and makes this span into a Hilbert
space. Let $k_w(z)=\frac{1}{z+w^*}$. Then for $h\in\mathbf
H_2(\mathbb C_r)$ and $f=h\circ\va$ we have:
\[
\begin{split}
\langle f, k_{\va(w)^*}(\va)\rangle&=\langle h,
k_{\va(w)^*}\rangle_{\mathbf H_2}\\
&=h(\va(w))\\
&=f(w).
\end{split}
\]
The result follows from the uniqueness of the reproducing kernel
Hilbert space associated to a given positive definite function.
\mbox{}\qed\mbox{}\\

\begin{Pn}
\begin{equation}
\label{sum123}
\mathbf H_2(\mathbb C_+)=\oplus_{n=1}^N  e_n(\mathbf H_2(\mathbb
C_+))(\va).
\end{equation}
\end{Pn}

{\bf Proof:} That the sum \eqref{sum123} is indeed orthogonal follows
from \eqref{qwerty}. Furthermore, let $f\in\mathbf H_2(\mathbb
C_+)$ be a finite linear span of kernels:
\[
f(z)=\sum_{j=1}^M\frac{a_j}{z+w_j^*}.
\]
Then, from \eqref{paris} we get
\[
f(z)=\sum_{n=1}^Ne_n(z)h_n(\va(z)),
\]
with
\[
h_n(z)=\sum_{j=1}^M\frac{a_j}{z+w_j^*}e_n(w_j)^*.
\]
We see that
\[
\begin{split}
\sum_{n=1}^N[h_n,h_n]&=\sum_{n=1}^N\sum_{j=1}^M
e_n(w_k)\frac{a_k^*a_j}{w_k+w_j^*}e_n(w_j)^*\\
&=[f,f].
\end{split}
\]
\mbox{}\qed\mbox{}\\

Note that neither $e_j$ nor $f(\va)$ belong to
$\mathbf H_2(\mathbb C_r)$. But we have:
\begin{Tm}
The maps
\[
C_jf(z)=e_j(z)f(\va(z))
\]
are continuous operators from the Hardy space $\mathbf
H_2(\mathbb C_r)$ into itself, and
\[
C_j^*\frac{1}{z+w^*}=\frac{e_j(w)^*}{z+\varphi(w)^*}.
\]
In particular
\[
\begin{split}
\sum_{n=1}^N C_jC_j^*&=I\\
C_k^*C_j&=\begin{cases}I\quad if\quad k=j\\
0\quad if \quad k\not = j.\end{cases}
\end{split}
\]
\end{Tm}

\section{Harmonic analysis of representations}
\setcounter{equation}{0}
\label{generalsetting}
As in the
introduction, we consider a function $k(z,w)$ positive definite
on a set $\mathcal U$, and a map $R$ from $\mathcal U$ into
itself. Recall that $R_n$ was defined by \eqref{eqrn}. We assume
that $k(z,w)$ is of the form
\begin{equation}
\label{sydney2012}
k(z,w)=1+t(z,w),
\end{equation}
where $t(z,w)$ is positive definite in $\mathcal U$. This is
equivalent to request that $\mathbb C$ is contractively included
in the reproducing kernel Hilbert spaces with reproducing kernel
$k(z,w)$. We set
\[
\Omega=\left\{z\in\mathcal U\,\,;\,\, \sum_{n=0}^\infty
|t(R_n(z),R_n(z))|<\infty\right\}.
\]
Note that this set may be empty, but that, in any case,
\[
R(\Omega)\subset\Omega.
\]

\begin{La}
Assuming that $\Omega\not =\emptyset$. Then the infinite product
\begin{equation}
\label{infinite1}
K(z,w)=\prod_{n=0}^\infty
(1+t(R_n(z),R_n(w)),\quad z,w\in\Omega,
\end{equation}
converges, and satisfies
\begin{equation}
\label{bastille}
K(z,w)=(1+t(z,w))K(R(z),R(w)),\quad z,w\in\Omega,
\end{equation}
\end{La}

{\bf Proof:} Since $t(z,w)$ is positive definite in $\Omega$ we have
\[
|t(R_n(z),R_n(w))|\le \sqrt{t(R_n(z),R_n(z))}\sqrt{t(R_n(w),R_n(w))},\quad z,w\in\Omega.
\]
The Cauchy-Schwarz inequality insures that
\[
\sum_{n=0}^\infty |t(R_n(z),R_n(w))|<\infty,
\]
and so the infinite product converges. Equation \eqref{bastille} follows from the definition of the
infinite product.
\mbox{}\qed\mbox{}\\

\begin{La}
Assume that $\Omega\not=\emptyset$. Then, $1\not\in\mathcal H(t)$
and
\begin{equation}
\label{inter}
\mathcal H(k)=\mathbb C\oplus\mathcal H(t).
\end{equation}
\end{La}

{\bf Proof:} By hypothesis, there exists $z\in\Omega_0$ such that
\begin{equation}
\label{unsw200712}
\lim_{n\rightarrow\infty} t(R_n(z),R_n(z))=0.
\end{equation}
Suppose that $1\in\mathcal H(t)$, and let $c=\|1\|^2_{\mathcal
H(t)}$. By formula \eqref{eqrk}, the kernel $t_1(z,w)$ defined by
\[
t(z,w)=\frac{1}{c}+t_1(z,w)
\]
is positive definite in $\Omega_0$. In particular we have
\[
t(R_n(z),R_n(z))\ge\frac{1}{c},\quad \forall n\in\mathbb
N\quad{\rm and}\quad \forall z\in\Omega_0,
\]
which contradicts \eqref{unsw200712}. From the decomposition
\eqref{sydney2012} we then have $1\in\mathcal H(k)$. Since
$\mathbb C\cap\mathcal H(t)=\left\{0\right\}$ we obtain
\eqref{inter}.
\mbox{}\qed\mbox{}\\

We now assume on $R$ the following two conditions: First,
\begin{equation}
\label{hypthR} \forall z\in\Omega,\quad n(z)\stackrel{\rm
def.}{=}{\rm Card}~\left\{\zeta\in\Omega,\,\:,\,
R(\zeta)=z\right\}<\infty,
\end{equation}
and one of the following two conditions:
\begin{equation}
\label{hypthR1} \forall z\in\Omega,\quad
\frac{1}{n(z)}\sum_{R(\zeta)=z}e_j(\zeta)\overline{e_k(\zeta)}=
\delta_{jk},\quad\forall j,k\in J,
\end{equation}

or

\begin{equation}
\label{hypthR12345} \forall z\in\Omega,\quad
\frac{1}{n(z)}\sum_{R(\zeta)=z}e_j(\zeta){e_k(\zeta)}=\delta_{jk},\quad\forall
j,k\in J,
\end{equation}
holds.

\begin{La}
Assume that \eqref{hypthR} is in force. Then:\\
$(a)$  If \eqref{hypthR1} is in force, the adjoint of the operator
$S_j$ is given by the formula
\begin{equation}
\label{sj*} (S_j^*f)(z)=\frac{1}{n(z)}
\sum_{\substack{\zeta\in\Omega_0\,\mbox{such}\\ \mbox{that} \,\,
R(\zeta)=z}} e_j(\zeta)^*f(\zeta).
\end{equation}
$(2)$ If \eqref{hypthR12345} is in force, the adjoint of the
operator $S_j$ is given by the formula
\begin{equation}
\label{sj**} (S_j^*f)(z)=\frac{1}{n(z)}
\sum_{\substack{\zeta\in\Omega\,\mbox{such}\\ \mbox{that} \,\,
R(\zeta)=z}} e_j(\zeta)f(\zeta).
\end{equation}

\end{La}

{\bf Proof:}\\
$(a)$ Using \eqref{bastille} we write for $z,w\in\Omega_0$
\[
\begin{split}
\frac{1}{n(z)}\sum_{R(\zeta)=z}(e_j(\zeta))^*
K(\zeta,w)&=\frac{1}{n(z)}\sum_{R(\zeta)=z}(e_j(\zeta))^*
(1+t(\zeta, w))K(R(\zeta),R(w))\\
&\hspace{-3cm} =\frac{1}{n(z)}\left(\sum_{k\in
J}\left(\sum_{R(\zeta)=z}
(e_j(\zeta))^*e_k(\zeta)\right)e_k(w)^*\right)K(R(\zeta),R(w))\\
&\hspace{-3cm}=(e_j(w))^*K(z,R(w)),\\
&\hspace{-3cm}=(S_j^*K_w)(z),
\end{split}
\]
by formula \eqref{mutl}, and where we have used \eqref{hypthR1} to go from the second to the third line.
Since the kernels are dense in $\mathcal H(K)$ and since $S_j^*$ is continuous,
the equality extends to all $f\in\mathcal H(K)$.

$(b)$ The proof is similar. One now has:
\[
\begin{split}
\frac{1}{n(z)}\sum_{R(\zeta)=z}(e_j(\zeta))
K(\zeta,w)&=\frac{1}{n(z)}\sum_{R(\zeta)=z}(e_j(\zeta))
(1+t(\zeta, w))K(R(\zeta),R(w))\\
&\hspace{-3cm} =\frac{1}{n(z)}\left(\sum_{k\in
J}\left(\sum_{R(\zeta)=z}
(e_j(\zeta))e_k(\zeta)\right)e_k(w)^*\right)K(R(\zeta),R(w))\\
&\hspace{-3cm}=(e_j(w))^*K(z,R(w)),\\
&\hspace{-3cm}=(S_j^*K_w)(z),
\end{split}
\]
\mbox{}\qed\mbox{}\\

An important case where the second set of conditions hold is
presented in \cite{aj_pnas}; see also Section \ref{sec7} below.

\begin{Tm}
\label{tmcuntz} Under hypothesis \eqref{hypthR} and
\eqref{hypthR1}, or \eqref{hypthR} and \eqref{hypthR12345} the
operators $(S_j)_{i\in J}$ satisfy the Cuntz relations.
\end{Tm}

{\bf Proof of Theorem \ref{tmcuntz}:} We first suppose that
\eqref{hypthR} and \eqref{hypthR1} hold. Let $i_0, j_0\in K$ and
$f\in\mathcal H(K)$. We have
\[
\begin{split}
(S_{i_0}^*S_{j_0}f)(z)&=\frac{1}{n(z)}
\sum_{\substack{\zeta\in\Omega\,\mbox{such}\\ \mbox{that} \,\,
R(\zeta)=z}}e_{i_0}(\zeta)^*
(S_{j_0}f)(\zeta)\\
&= \frac{1}{n(z)}\sum_{\substack{\zeta\in\Omega\,
\mbox{such}\\ \mbox{that} \,\, R(\zeta)=z}}e_{i_0}(\zeta)^*e_{j_0}(\zeta)f(R(\zeta))\\
&=
\left(\frac{1}{n(z)}\sum_{\substack{\zeta\in\Omega\,\mbox{such}\\
\mbox{that} \,\, R(\zeta)=z}}e_{i_0}(\zeta)^*e_{j_0}(\zeta)\right)
f(z)\\
&=\delta_{i_0,j_0}f(z),\quad\mbox{\rm thanks to \eqref{hypthR1}}.
\end{split}
\]
We now assume that \eqref{hypthR} and \eqref{hypthR12345} are in
force. Then,
\[
\begin{split}
(S_{i_0}^*S_{j_0}f)(z)&=\frac{1}{n(z)}
\sum_{\substack{\zeta\in\Omega\,\mbox{such}\\ \mbox{that} \,\,
R(\zeta)=z}}e_{i_0}(\zeta)
(S_{j_0}f)(\zeta)\\
&= \frac{1}{n(z)}\sum_{\substack{\zeta\in\Omega\,
\mbox{such}\\ \mbox{that} \,\, R(\zeta)=z}}e_{i_0}(\zeta)e_{j_0}(\zeta)f(R(\zeta))\\
&=
\left(\frac{1}{n(z)}\sum_{\substack{\zeta\in\Omega_0\,\mbox{such}\\
\mbox{that} \,\, R(\zeta)=z}}e_{i_0}(\zeta)e_{j_0}(\zeta)\right)
f(z)\\
&=\delta_{i_0,j_0}f(z),\quad\mbox{\rm thanks to
\eqref{hypthR12345}}.
\end{split}
\]
\mbox{}\qed\mbox{}\\

\section{An orthogonal basis}
\setcounter{equation}{0} In this section we show that for anyone
of the representations of a fixed$\mathcal O_n$  in some Hilbert
space
$\mathcal H$, one may naturally construct an associated orthonormal basis
(ONB) in $\mathcal H$. We will explore its implications for the analysis
of kernel Hilbert spaces with special view to those arising from
the iterated function systems in Julia set theory.\\

From the infinite product representation \eqref{infinite1} of $K(z,w)$ we see that
\[
K(z,w)=1+K_1(z,w),
\]
where $K_1(z,w)$ is positive definite in $\Omega_0$. Furthermore,
the function $\mathbf 1$:
\[
\mathbf 1(z)\equiv 1, \quad z\in\Omega,
\]
belongs to $\mathcal H(K)$, and in particular
\[
e_j=S_j(\mathbf 1)\in\mathcal H(K),\quad\forall j\in J.
\]
In this section we wish to express $K$ in the representation of
the form \eqref{eqrk} for an appropriate basis expressed in terms
of the function $\mathbf 1$ and of the $S_j$. We set $N={\rm
dim}~\mathcal H(k)$, that is the cardinal of $J$ (possibly, $N\ge
\aleph_0$), and consider $V$ the tree with at each vertex $N$
edges associated to $(S_j)_{j\in J}$. On the vertices of the tree
we have the functions
\begin{equation}
\label{eq:bv}
b_v(z)=\left(S_{i_0}S_{i_1}\cdots S_{i_N}\mathbf 1
\right)(z),
\end{equation}
where $N=0,1,2,\ldots$ and the $i_j$ belong to the index set $J$,
formed from an iterated application of the $S_{i_j}$.

\begin{Tm}
The functions $(b_v)_{v\in V}$ form an orthonormal basis of $\mathcal H(K)$ and
it holds that:
\begin{equation}
\label{mabillon}
K(z,w)=\sum_{v\in V} b_v(z)b_v(w)^*
\end{equation}
where  $b_v$ is given by \eqref{eq:bv}.
\end{Tm}

{\bf Proof:} The Cuntz relations readily imply that the $(b_v)_{v\in V}$ form an orthonormal system. We need to see that it is complete. To see this it is enough to check directly that \eqref{mabillon} holds. Let \eqref{formule_k} be a
representation of $k$ in terms of an orthonormal basis $(e_j)_{j\in J}$ of $\mathcal H(k)$.
The infinite product \eqref{infinite1} is equal to a sum of elements of the form $f(z)f(w)^*$, where $f$ is of the form
\[
f(z)=e_{i_1}(z)e_{i_2}(R(z))e_{i_3}(R^2(z))\cdots e_{i_M}(R^{M-1}(z)),
\]
where $M=1,2,\ldots$ and the $i_j$ belong to the index set $J$.
Indeed, \eqref{infinite1}
is equal to the limit
\[
K(z,w)=\lim_{N\longrightarrow\infty}\prod_{n=0}^N
\left(\sum_{j\in J}e_j(R^n(z))(e_j(R^n(w))^*)\right),\quad z,w\in\Omega.
\]
For a given $N$ we have
\[
\begin{split}
\prod_{n=0}^N \left(\sum_{j\in J}e_j(R^n(z)(e_j(R^n(w))^*\right)&=\\
&\hspace{-3cm}=\sum_{(i_1,\ldots, i_N)\in J^N}
e_{i_1}(z)e_{i_2}(R(z))e_{i_3}(R^2(z))\cdots e_{i_N}(R^{N}(z)),
\end{split}
\]
that is
\[
\prod_{n=0}^N \left(\sum_{j\in
J}e_j(R^n(z)(e_j(R^n(w))^*\right)=\sum_{|v|=N+1} b_v(z)b_v(w)^*,
\]
where we have denoted by $|v|$ the length of the path $v$
starting at the origin. The result follows since the infinite
product converges. Finally, by definition of the operators
$S_{i_j}$ we have:
\[
\begin{split}
e_{i_0}(z)e_{i_1}(R(z))e_{i_2}(R^2(z))\cdots e_{i_N}(R^{N}(z))
&=\\
&\hspace{-1.5cm}=
S_{i_0}\left(e_{i_1}(\cdot)e_{i_2}(R(\cdot))\cdots e_{i_N}(R^{N-1}(\cdot))\right)(z)\\
&\hspace{-1.5cm}= S_{i_0}\left(S_{i_2}\left(e_{i_3}(\cdot)\cdots
e_{i_{N}}(R^{N-2}
(\cdot))\right)\right)(z)\\
&\hspace{2mm}\vdots\\
& \hspace{-1.5cm}= \left(S_{i_0}S_{i_1}S_{i_2}\cdots
S_{i_N}\right)(\mathbf 1).
\end{split}
\]
This concludes the proof.
\mbox{}\qed\mbox{}\\

\section{Example: A Julia set}
\setcounter{equation}{0}
\label{sec7}
We consider $P(z)=z^2-1$ and
\[
R(z)=P(P(z))=z^4-2z^2.
\]
We check below that the conditions \eqref{hypthR} and
\eqref{hypthR12345} are in force. We first define
\[
\Omega=\left\{w\in\mathbb C\,\,\mbox{\rm such that}\,\,
(R_n(w))_{n\in\mathbb N_0}\in \ell_2\right\}.
\]
For $z,w\in\Omega$ we set
\begin{equation}
K(z,w)=\prod_{n=0}^\infty (1+R_n(z)R_n(w)^*).
\label{kernelR}
\end{equation}
\begin{Pn}
The infinite product \eqref{kernelR} converges in $\Omega$ to a
function $K(z,w)$ which is positive definite there. Furthermore,
$K$ satisfies the equation
\[
K(z,w)=E(z,w)K_R(z,w),\quad z,w\in\Omega,
\]
with
\begin{equation}
\label{eqrecurs} K_R(z,w)=K(R(z),R(w))\quad and\quad
E(z,w)=1+zw^*.
\end{equation}
\end{Pn}

{\bf Proof:} Since $\ell_2\subset\ell_1$ the Cauchy-Schwarz
inequality insures that
\[
\sum_{n=0}^\infty |R_n(z)R_n(w)^*|<\infty,\quad z,w\in\Omega,
\]
and so the infinite product converges there. The limit is
positive definite in $\Omega$ since each of the factor is
positive definite there and since a convergent product of
positive definite functions is positive definite. Finally,
equation \eqref{eqrecurs} is clear from the infinite product
representation of $K$.
\mbox{}\qed\mbox{}\\

\begin{Pn}
$\Omega$ is equal to the Fatou set at $0$.
\end{Pn}

{\bf Proof:} One direction is clear. If $z\in\Omega$, then
$\lim_{n\rightarrow\infty} R_n(z)=0$, and so $z$ is in the Fatou
set. Conversely, let $z$ be in the Fatou set. Then there is $n_0$
such that
\[
n\ge n_0\quad\longrightarrow |R_n(z)|<\frac{1}{2}.
\]
But $R_{n+1}(z)=(R_n(z))^2((R_n(z))^2-2)$, and so
\[
|R_{n+1}(z)|=\le\frac{3}{4}|R_n(z)|,
\]
and $(R_n(z))_{n\in\mathbb N}\in\ell_2$.
\mbox{}\qed\mbox{}\\
\begin{Pn}
Let $z\in\Omega$. The equation $R(\zeta)=z$ has four solutions in
$\Omega$ \label{pn1}
\end{Pn}

{\bf Proof:}

$R(\zeta)=z$ reads
\begin{equation}
\label{eqzeta}
\zeta^4-2\zeta^2-z=0,
\end{equation}
and so hypothesis \eqref{hypthR}
\[
\forall z\in\mathbb C,\quad n(z)\stackrel{\rm def.}{=}{\rm
Card}~\left\{\zeta\in\Omega_0,\,\:,\, R(\zeta)=z\right\}<\infty
\]
holds with $n(z)=4$  by the fundamental theorem of algebra. That
the solutions belong to $\Omega$ follows from the fact that
$P^{-1}(\Omega)=\Omega$ (the inverse image of the Fatou set is
the Fatou set; see \cite{MR1128089}, \cite{MR1721240}), and so $R^{-1}(\Omega)=\Omega$.
\mbox{}\qed\mbox{}\\

We note that the sums \eqref{hypthR12345} now read, with
$e_1(z)=1$ and $e_2(z)=z$
\begin{equation}
\label{juliarel}
\begin{split}
\sum_{R(\zeta)=z}e_j(\zeta)&=0,\\
\sum_{R(\zeta)=z}e_j^2(\zeta)&=4,\\
\sum_{R(\zeta)=z, k\not=j}e_k(\zeta)e_j(\zeta)&=0.
\end{split}
\end{equation}
The first one is in force because the coefficient of $\zeta$ is
$0$ in \eqref{eqzeta}. The third one reduces to the first one
since $e_1(z)=1$. To verify the second one, let $x(z)$ be a
complex number such that $x(z)^2=1+z$. Then
\[
\zeta^2=1\pm x(z),
\]
and the second equation follows.\\

\begin{equation}
\begin{split}
S_0f(z)&=f(R(z))\\
S_1(z)&=zf(R(z)).
\end{split}
\end{equation}
\begin{Pn}
$S_0$ and $S_1$ are bounded operators from $\mathcal H(K)$ into
itself. They satisfy
\begin{equation}
S_0S_0^*+S_1S_1^*=I_{\mathcal H(K)}.
\end{equation}
\end{Pn}

\begin{La}
We have
\begin{equation}
\label{s0*}
(S_0^*f)(z)=\frac{1}{n(z)}\sum_{\substack{\zeta\in\Omega\,\mbox{such}\\
\mbox{that} \,\, R(\zeta)=z}}f(\zeta),
\end{equation}
and
\begin{equation}
\label{s1*}
(S_1^*f)(z)=\frac{1}{n(z)}\sum_{\substack{\zeta\in\Omega\,
\mbox{such}\\ \mbox{that} \,\, R(\zeta)=z}}{\zeta}f(\zeta)
\end{equation}
\end{La}

{\bf Proof:} We follow the argument in \cite{aj_pnas}. To prove
\eqref{s0*} we write:
\[
\begin{split}
\frac{1}{n(z)}\sum_{R(\zeta)=z}K(\zeta,w)&=\frac{1}{n(z)}\sum_{R(\zeta)=z}
(1+\zeta w^*)K(R(\zeta),R(w))\\
&=\left(1+\frac{\left(\sum_{R(\zeta)=z}\zeta\right) w^*}{n(z)}K(z,R(w))\right)\\
&=K(z,R(w)),\quad\mbox{\rm since \eqref{hypthR} is in force}\\
&=(S_0^*K_w)(z).
\end{split}
\]
The result follows by density since $S_0^*$ is continuous. The
argument for $S_1^*$ is as follows:
\[
\begin{split}
\frac{1}{n(z)}\sum_{R(\zeta)=z}\zeta K(\zeta,w)&=\frac{1}{n(z)}
\sum_{R(\zeta)=z}\zeta(1+\zeta w^*)K(R(\zeta),R(w))\\
&=\left(1+\frac{\left(\sum_{R(\zeta)=z}\zeta^2\right) w^*}{n(z)}K(z,R(w))\right)\\
&=w^*K(z,R(w)),\quad\mbox{\rm since \eqref{hypthR12345} is in force}\\
&=(S_1^*K_w)(z).
\end{split}
\]
\mbox{}\qed\mbox{}\\

\begin{Tm}
Assume \eqref{hypthR} and \eqref{hypthR1} in force. Then the pair
of operators $(S_0,S_1)$ satisfies the Cuntz relations in
$\mathcal H(K)$.
\end{Tm}

{\bf Proof:} We have
\[
(S_0^*S_0f)(z)=\frac{1}{n(z)}\sum_{R(\zeta)=z}(S_0f)(z)=
\frac{1}{n(z)}\sum_{R(\zeta)=z}f(R(\zeta))=f(z),
\]
and
\[
\begin{split}
(S_0^*S_1f)(z)&=\frac{1}{n(z)}\sum_{R(\zeta)=z}(S_1f)(\zeta)\\
&=\frac{1}{n(z)}\sum_{R(\zeta)=z}\zeta f(R(\zeta))\\
&=\frac{1}{n(z)}\left(\sum_{R(\zeta)=z}\zeta\right)f(z)=0.
\end{split}
\]
Finally, the computation for $S_1^*S_1$ is as follows:.
\[
\begin{split}
(S_1^*S_1f)(z)&=\frac{1}{n(z)}\sum_{R(\zeta)=z}\zeta(S_1f)(\zeta)\\
&=\frac{1}{n(z)}\sum_{R(\zeta)=z}\zeta^2 f(R(\zeta))\\
&=\frac{1}{n(z)}\sum_{R(\zeta)=z}\zeta^2 f(z)\\
&=f(z),
\end{split}
\]
where we have used the second equality in \eqref{juliarel}.
\mbox{}\qed\mbox{}

\bibliographystyle{plain}
\def\cprime{$'$} \def\lfhook#1{\setbox0=\hbox{#1}{\ooalign{\hidewidth
  \lower1.5ex\hbox{'}\hidewidth\crcr\unhbox0}}} \def\cprime{$'$}
  \def\cfgrv#1{\ifmmode\setbox7\hbox{$\accent"5E#1$}\else
  \setbox7\hbox{\accent"5E#1}\penalty 10000\relax\fi\raise 1\ht7
  \hbox{\lower1.05ex\hbox to 1\wd7{\hss\accent"12\hss}}\penalty 10000
  \hskip-1\wd7\penalty 10000\box7} \def\cprime{$'$} \def\cprime{$'$}
  \def\cprime{$'$} \def\cprime{$'$}

\end{document}